\documentclass[12pt]{article}
\usepackage[reqno]{amsmath}
\usepackage{amssymb, theorem, enumerate}
\usepackage{graphicx}

\theoremstyle{change}
{\theorembodyfont{\slshape}
\newtheorem{theorem}{Theorem.}[section]
\newtheorem{lemma}[theorem]{Lemma.}
\newtheorem{corollary}[theorem]{Corollary.}}
\theorembodyfont{\rmfamily}

\newcommand\lref[1]{Lemma~\ref{lem:#1}}

\newcommand\cref[1]{Corollary~\ref{cor:#1}}
\newcommand\sref[1]{Section~\ref{sec:#1}}

\def\proof{\noindent{{\sl Proof. }}}
\def\sqr#1#2{{\vbox{\hrule height.#2pt
    \hbox{\vrule width.#2pt height#1pt \kern#1pt
        \vrule width.#2pt}\hrule height.#2pt}}}
\def\eqed{\sqr53}
\def\qed{%
    \ifmmode\eqno\eqed
    \else\nobreak\ \hfill\eqed\medbreak\fi}

\newcommand\be{\beta}

\newcommand\De{\Delta}

\newcommand\ga{\gamma}

\newcommand\sg{\sigma}

\newcommand\cC{{\mathcal C}}

\newcommand\ints{{\mathbb Z}}

\newcommand\comp[1]{{\mkern2mu\overline{\mkern-2mu#1}}}
\newcommand\diff{\mathbin{\mkern-1.5mu\setminus\mkern-1.5mu}}

\newcommand\sbs{\subseteq}

\DeclareMathOperator{\supp}{supp}

\newcommand\pmat[1]{\begin{pmatrix} #1 \end{pmatrix}}

\newcommand\one{{\bf1}}

\DeclareMathOperator{\tr}{tr}

\DeclareMathOperator\wt{wt}

\newcommand\dv{D}

\title{Perfect State Transfer in Cubelike Graphs}
\author{Wang-Chi Cheung and Chris Godsil\\[2pt]
Combinatorics \& Optimization\\
University of Waterloo}

\date{\small{AMS Classification: 05C50, 06E99, 81P68 \\
Keywords: binary codes, cubelike graph, perfect state transfer}}

\begin{document}
\maketitle

\begin{abstract}
	Suppose $C$ is a subset of non-zero vectors from the vector space 
	$\ints_2^d$. The \textsl{cubelike graph} $X(C)$ has $\ints_2^d$ as 
	its vertex set, and two elements of $\ints_2^d$ are adjacent if their 
	difference is in $C$.
	If $M$ is the $d\times |C|$ matrix with the elements of $C$ as its columns,
	we call the row space of $M$ the \textsl{code} of $X$.
	We use this code to study perfect state transfer on cubelike graphs.  
	Bernasconi et al have shown that perfect state transfer occurs on $X(C)$
	at time $\pi/2$ if and only if the sum of the elements of $C$ is not zero.  
	Here we consider what happens when this sum is zero.
	We prove that if perfect state transfer occurs
	on a cubelike graph, then it must take place at time $\tau=\pi/2\dv$, 
	where $\dv$ is the greatest common divisor of the 
	weights of the code words. We show that perfect state transfer occurs at time
	$\pi/4$ if and only if $\dv=2$ and the code is self-orthogonal.
\end{abstract}

\section{Introduction}

Let $X$ be a graph on $v$ vertices with adjacency matrix $A$.  We define a 
transition operator $H(t)$ by
\[
	H(t) :=\exp(iAt)
\]
This operator is unitary and in quantum computing it determines a continuous
quantum walk \cite{Kempe:2003p4062}.  We say that \textsl{perfect state transfer} 
from vertex $u$ to vertex $v$ occurs at time $\tau$ if $u\ne v$ and
\[
	|H(\tau)_{u,v}| =1
\]
If we have
\[
	|H(\tau)_{u,u}| =1
\]
then we say that $X$ is \textsl{periodic at $u$} with period $\tau$. We say $X$
is \textsl{periodic} with period $\tau$ if it is periodic at each vertex with
period $\tau$. Of these two concepts, perfect state transfer is the one of
physical interest but we will see that periodicity is closely related to it.

Perfect state transfer was studied in detail by 
Christandl et al \cite{Christandl:2005p4044} 
where they showed that, in the $d$-cube, perfect state transfer occurs at 
time $\pi/2$ from each vertex to the unique vertex at distance $d$ from it.
For a recent survey on perfect state transfer see \cite{KendTam2011}.

The $d$-cube is an example of a Cayley graph for $\ints_2^d$.
A Cayley graph $X(C)$ for $\ints_2^d$ has the binary vectors of
length $d$ as its vertices, with two vertices adjacent if and only if their
difference lies in some specified subset $C$ of $\ints_2^d\diff\{0\}$. (The set $C$
is the \textsl{connection set} of the Cayley graph.) If we choose $C$ to consist
of the $d$ vectors from the standard basis of $\ints_2^d$, then the cubelike
graph $X(C)$ is the $d$-cube. In \cite{Facer:2008p4037} Facer, Twamley and Cresser 
showed that perfect state transfer occurs in a special class of
Cayley graphs for $\ints_2^d$ that includes the $d$-cube, and this was extended 
to an even larger class of graphs in \cite{Bernasconi:2008p4036} by Bernasconi, 
Godsil and Severini.  

If we let $\sg$ denote the sum of the elements of $C$, then the main result 
of \cite{Bernasconi:2008p4036} is that, if $\sg\ne0$, then at time $\pi/2$
we have perfect state transfer from $u$ to $u+\sg$, for each vertex $u$.
Our goal in this paper is to study the situation when $\sg=0$; we find a surprising
connection to binary codes.

\section{Perfect State Transfer}

If $u\in\ints_2^d$, then the map
\[
	x\mapsto x+u
\]
is a permutation of the elements of $\ints_2^d$, and thus it can be represented
by a $2^d\times 2^d$ permutation matrix $P_u$.  We note that $P_0=I$,
\[
	P_u P_v =P_{u+v}
\]
and so $P_u^2=I$.  We also see that $\tr(P_u)=0$ if $u\ne0$ and
\[
	\sum_{u\in\ints_2^d}P_u =J.
\]

\begin{lemma}\label{lem:adj-perms}
	If $C\sbs\ints_2^d\diff0$ and $X$ is the cubelike graph with connection set
	$C$, then $A(X)=\sum_{u\in C}P_u$.\qed
\end{lemma}

If $\sg$ is the sum of the elements of $C$, then
\[
	P_\sg =\prod_{u\in C}P_u.
\]

\begin{lemma}\label{lem:Hprod}
	If $H(t)$ is the transition operator of the cubelike graph $X(C)$, then
	$H(t)=\prod_{u\in C}\exp(itP_u)$.
\end{lemma}

\proof
If matrices $M$ and $N$ commute then
\[
	\exp(M+N) =\exp(M)\exp(N)
\]
Since $A=\sum_{u\in C}P_u$ and since the matrices $P_u$ commute, the lemma follows.\qed

Suppose $P$ is a matrix such that $P^2=I$.  Then
\[
	\exp(itP) =I +itP -\frac{t^2}{2!}I -i\frac{t^3}{3!}P +\frac{t^4}{4!}I +\cdots
\]
and hence
\[
	\exp(itP) = \cos(t)I +i\sin(t)P.
\]
If $P$ is a permutation matrix we see that
\[
	\exp(\pi iP) = -I,\quad \exp\Bigl(\frac12\pi iP\Bigr) = iP.
\]
This implies that we have perfect state transfer on $K_2$ at time $\pi/2$, 
and that $K_2$ is periodic with period $\pi$.

If $H$ is the transition operator for a Cayley graph of an abelian group
then the argument used above shows that $H$ can be factorized
as a product of transition operators for a collection of perfect matchings
and 2-regular subgraphs. Unfortunately this does not seem to allow us to
derive usual information about state transfer.

\medbreak
Now we present a new and very simple proof of Theorem~1 from Bernasconi et al
\cite{Bernasconi:2008p4036}.

\begin{theorem}\label{thm:abgs}
	Let $C$ be a subset of $\ints_2^d$ and let $\sg$ be the sum of the elements
	of $C$.  If $\sg\ne0$, then perfect state transfer occurs in $X(C)$ from
	$u$ to $u+\sg$ at time $\pi/2$.  If $\sg=0$, then $X$ is periodic with 
	period $\pi/2$.
\end{theorem}

\proof
Let $H(t)$ be the transition operator associated with $A$.  Then by \lref{Hprod}
\[
	H(t) =\prod_{u\in C}\exp(itP_u).
\]
From our remarks above
\[
	\exp(itP_u) =\cos(t)I +i\sin(t)P_u
\]
and therefore
\[
	H(\pi/2) =\prod_{u\in C} iP_u =i^{|C|} P_\sg.
\]
This proves both claims.\qed

This result is very natural, but clearly raises the question of whether
we can have perfect state transfer when $\sg\ne0$. We will see that we can.

We show how to use these ideas to arrange for perfect state transfer from $0$ to a 
specified vertex $u$ in a cubelike graph. Assume we have cubelike graph
with connection set $C$ and let $\sg$
be the sum of the elements of $C$. If $\sg=u$ then we already have transfer to $u$.
First assume $\sg=0$. If $u\in C$ let $C'$ denote $C\diff u$; if $u\notin C$ let
$C'$ be $C\cup u$. In both cases the sum of the elements of $C'$ is $u$ and we're done.
If $\sg\ne 0$, replace $C$ by $(C\diff\sg)$, now we are back in the first case.
We can summarize this as follows. Let $S\oplus T$ denote the symmetric difference of sets
$S$ and $T$.

\begin{lemma}\label{lem:}
    If $u$ is a vertex in the cubelike graph $X(C)$, then there is a connection
    set $C'$ such that $|C\oplus C'|\le2$ and we have perfect state transfer from
    $0$ to $u$ in $X(C')$ at time $\pi/2$.\qed
\end{lemma}

\section{The Minimum Period}
\label{sec:period}

In this section we determine the minimum period of a cubelike graph.

We consider the spectral decomposition of the adjacency matrix of a cubelike 
graphs.  If $a\in\ints_2^d$, then the function
\[
	x\mapsto (-1)^{a^Tx}
\]
is both a character of $\ints_2^d$ and an eigenvector of $X(C)$ with eigenvalue
\[
	\sum_{c\in C}(-1)^{a^Tc}.
\]
Let $M$ be the matrix with the elements of $C$ as its columns.  Its row space is
a binary code, and if $\wt(x)$ denotes the Hamming weight of $x$, the above 
eigenvalue is equal to
\[
	|C| - 2\wt(a^TM).
\]
Thus the weight distribution of the code determines the eigenvalues of $X(C)$,
and also their multiplicities.

As a pertinent example we offer
\[
    M=\pmat{
    0 & 0 & 0 & 0 & 0 & 0 & 0 & 1 & 1 & 1 & 1 \\
    0 & 0 & 0 & 0 & 0 & 1 & 1 & 0 & 0 & 1 & 1 \\
    0 & 0 & 0 & 1 & 1 & 0 & 0 & 0 & 0 & 1 & 1 \\
    0 & 1 & 1 & 0 & 0 & 0 & 0 & 0 & 0 & 1 & 1 \\
    1 & 0 & 1 & 0 & 1 & 0 & 1 & 0 & 1 & 0 & 1}
\]
which has weight enumerator
\[
    x^{11} + 10x^7y^4 + 16x^5y^6 + 5x^3y^8,
\]
from which we learn that the weights of its code words are $0$, $4$, $6$ and $8$.
The eigenvalues of the associated cubelike graph are
\[
    11,\ 3,\ -1,\ -5
\]
with respective multiplicities
\[
    1,\ 10,\ 16,\ 5.
\]

If we define the $2^d\times 2^d$ matrix $E_a$ by
\[
	(E_a)_{u,v} := 2^{-d} (-1)^{a^T(u+v)}
\]
then $E_a^2=E_a$ and, if $a\ne b$, then $E_aE_b=0$.  The columns of $E_a$ are
eigenvectors for $X(C)$ with eigenvalue $|C|-2\wt(a^TM)$, and if $m=|C|$ we have
\[
	A =\sum_a (m-2\wt(a^TM))E_a.
\]
More significantly
\begin{align*}
	\exp(iAt) &=\sum_a \exp(i(m-2\wt(a^TM))t)E_a\\
			&=\exp(imt)\sum_a \exp(-2it\wt(a^TM))E_a
\end{align*}

\begin{lemma}\label{lem:period}
	Let $X$ be a cubelike graph and let $\dv$ be the greatest common divisor
	of the weights of the words in its code. Then the minimum period of $X$ 
	is $\pi/\dv$.
\end{lemma}

\proof
If $X$ is periodic with period $\tau$, from Theorem~4.1 in \cite{Godsil:2008p4050} 
we know there is a complex scalar $\ga$
with norm 1 such that $H(\tau)_{u,u}=\ga$, for any vertex $u$.  Therefore
\[
	\ga = 2^{-d}\exp(mi\tau) \sum_a \exp(-2i\tau\wt(a^TM)).
\] 
This shows that $\ga\exp(-mi\tau)$ is the average of the $2^d$ terms in the 
above sum. Since each of these terms lies on the complex unit circle and, 
since $\ga\exp(-m\tau)$ lies
on the unit circle, we conclude that for all choices of $a$,
\[
	\ga\exp(-mi\tau) =\exp(-2i\tau\wt(a^TM))
\]
Set $q$ equal to $\exp(-2i\tau)$.  Then for any $a$ and $b$ we have
\[
	q^{\wt(a^TM)} =q^{\wt(b^TM)}
\]
and accordingly
\[
	\Bigl(q^{\dv}\Bigr)^{(\wt(b^TM) -\wt(a^TM))/\dv} =1.
\]
Since the set of nonzero integers of the form $(\wt(b^TM) -\wt(a^TM))/\dv$ 
is coprime, we conclude that $q^\dv=1$.  Hence $\exp(-2i\dv\tau)=1$ and therefore
\[
	\tau =\frac{\pi}{\dv}.\qed
\]

It follows from our calculations that $\ga=\exp(im\pi/\dv)$.

\section{Characterizing State Transfer}

\begin{theorem}\label{thm:pst}
    Let $X$ be a cubelike graph with matrix $M$ and suppose $u$ is a vertex
    in $X$ distinct from $0$. Then the following are equivalent:
    \begin{enumerate}[(a)]
        \item 
        There is perfect state transfer from $0$ to $u$ at time $\pi/2\De$.
        \item
        All words in $C$ have weight divisible by $\De$ and
        $\De^{-1}\wt(a^TM)$ and $a^Tu$ have the same parity for all vectors $a$.
        \item
        $\De$ divides $|\supp(u)\cap\supp(v)|$ for any two code words $u$ and $v$.
    \end{enumerate}
\end{theorem}

\proof 
We start by proving that (a) and (b) are equivalent. Perfect state transfer
occurs at time $\pi/2\De$ if and only if there is a complex
scalar $\be$ of norm 1 and a permutation matrix $T$ of order two and with trace
zero such that
\[
	H(\pi/2\De) =\be T.
\]
Now
\[
	(H(\pi/2\De))_{0,u} =\exp(im\pi/2\De) \sum_a    
	    \exp(-i\pi\wt(a^TM)/\De)(E_a)_{0,u}
\]
and
\[
	(E_a)_{0,u} = 2^{-d}(-1)^{a^Tu},
\]
consequently
\begin{align*}
	\be\exp(-im\pi/\De) &= 2^{-d}\sum_a \exp(-i\pi\wt(a^TM)/\De) (-1)^{a^Tu}\\
		&=2^{-d}\sum_a (-1)^{\wt(a^TM)/\De}(-1)^{a^Tu}.
\end{align*}
Here the left side of this equation has absolute value 1 and the right side is the 
average of $2^d$ numbers of absolute value 1, so the left side is $\pm1$ and the summands
on the right all have the same sign. So this equation holds if and only if, for all $a$ 
we have
\[
	\frac{\wt(a^TM)}{\De} = a^Tu,\quad (\mathrm{modulo}\ 2).
\]
Now this holds if and only if, modulo 2,
\[
	\frac{\wt((a+b)^TM)}{\De} = \frac{\wt(a^TM)}{\De}+\frac{\wt(b^TM)}{\De}
\]
for all $a$ and $b$.  This holds in turn if and only if, for any
two code words $u$ and $v$, we have that
\[
	\wt(u+v) = \wt(u)+\wt(v)\quad \textrm{(mod $2\De$)}
\]
and this holds if and only if
\[
	|\supp(u)\cap\supp(v)| = 0 \quad\textrm{(mod $\De$)}.\qed.
\]

Suppose $\cC$ is a binary code with generator matrix $M$. Let $M'$ denote $M$ viewed as 
a $01$-matrix over $\ints$ and let $\De$ be the gcd of the entries in $M'\one$. 
Then the entries of $\De^{-1}M'\one$ are integers, not all even, and we define the 
image of this vector in $\ints_2$ to be the \textsl{center} of $\cC$. Note that $\De$ 
is the gcd of the weights of the code words formed by the rows of $M$ and, if
$\De$ is odd, then the centre of $\cC$ is equal to $M\one$.

\begin{corollary}
    Suppose $X$ is a cubelike graph and $u$ is a vertex in $X$ distinct from 0.
    If we have perfect state transfer from $0$ to $u$ at time $\pi/2\De$,
    then $\De$ is the divisor of the code of $X$, and $u$ is its centre.\qed
\end{corollary}

\proof
Clearly $\dv|\De$. Since the size of intersection of the supports of two code words is
divisible by $\De$, it follows by induction that the weight of any sum of $k$
rows of $M$ has weight divisible by $\De$ and hence $\De|\dv$.\qed

Suppose $x$ and $y$ are binary vectors and $\De$ divides the weight of $x$, $y$
and $x+y$. If
\[
    \wt(x)=a+b,\quad \wt(y)=a+c;\quad \wt(x+y) = b+c
\]
then, modulo $\De$,
\begin{align*}
    a+b\phantom{+c} &= 0\\
    a\phantom{+b}+c &= 0\\
    \phantom{a+}b+c &= 0.
\end{align*}
This implies that, modulo $\De$,
\[
    2a =2b =2c =0.
\]
It follows that the odd integer $d$ divides the weight of each word in
a binary code if and only if, for any two words $x$ and $y$, the size of
$\supp(x)\cap\supp(y)$ is divisible by $d$.

\section{Examples}

A code is \textsl{even} if $\dv$ is even and \textsl{doubly even} if
$\dv$ is divisible by four. If $C$ is even and the size of the intersection
of any two codes is even, then $C$ is self-orthogonal. 
Note that since our graphs are simple, their generator matrices cannot
have repeated columns. (Using the standard terminology our codes
are projective or, equivalently, the minimum distance of the dual is at 
least three.) So cubelike
graphs with perfect state transfer at time $\pi/4$ correspond
to self-orthogonal projective binary codes that are even but not doubly even.

Unpublished computations by Gordon Royle have provided a complete list
of the cubelike graphs on 32 vertices. Analysis of the graphs in this list that 
show there are exactly six cubelike graphs on 32 vertices for which the codes are 
self-orthogonal and even but not doubly-even. The example in \sref{period} is the one
of these with least valency. These six graphs split into three pairs, 
each the complement of the other. In general, if perfect state transfer occurs on a
graph it need not occur on its complement. In our case it must, as the following 
indicates. We use $\comp{X}$ to denote the complement of $X$.

\begin{lemma}
   If $X$ is a cubelike graph with at least eight vertices then prefect state
   transfer occurs on $X$ if and only if it occurs on $\comp{X}$. 
\end{lemma}

Since $A(\comp{X})=J-I-A(X)$ we have.
\[
    H_{\comp{X}}(t) =\exp(it(J-I-A)).
\] 
If $X$ is regular then $J$ and $A$ commute and
\[
    H_{\comp{X}}(t) = \exp(it(J-I))H_X(-t)
\]
and hence
\[
    H_{\comp{X}}(\pi/k) = \exp(-i\pi/k) \exp((\pi i/k)J) H_X(\pi/k)^{-1}
\]
If $|V(X)n=$ then the eigenvalues of $J$ are 0 and $n$ and $\exp((\pi/k)J)=I$
provided $n/k$ is even.

There are a further six cubelike graphs on 32 vertices whose codes
are doubly even. A doubly even code is necessarily self-orthogonal.
If perfect state transfer occurs at time $\tau$, then Lemma~5.2 in 
\cite{Godsil:2008p4050}
yields that $\tr(H(\pi/4))$ must be zero, and using this we can show that 
perfect state transfer does not occur on these graphs. Thus we do not have 
examples of cubelike graphs with $\dv>2$ where perfect state transfer occurs.

If $M$ and $N$ are binary matrices, their \textsl{direct sum} is the matrix
\[
    \pmat{M&0\\ 0&N}
\]
and the code of this matrix is the direct sum of the codes of $M$ and $N$.
If the code of $M$ is self-orthogonal and even but not doubly even, then
the direct sum of two copies of this code is all these things too.
If $X$ and $Y$ are the cubelike graphs belonging to $M$ and $N$, then 
the cubelike graph belonging to the direct sum of $M$ and $N$ is the
Cartesian product of $X$ and $Y$. The transition matrix of the Cartesian product
of $X$ and $Y$ is $H_X\otimes H_Y$. One consequence is that we do have infinitely
many examples of cubelike graphs admitting perfect state transfer at time $\pi/4$.
We would very much like to know if perfect state transfer on cubelike graphs
could occur with a period less than $\pi/4$.

\end{document}